\documentclass[final,12pt,3p]{elsarticle}

\usepackage{graphicx}
\usepackage{amssymb}
\usepackage{amsthm}
\usepackage{amsmath}
\usepackage{algorithm}
\usepackage{algorithmic}
\usepackage{tikz}
\usepackage{lineno}
\usepackage{hyperref}
\usepackage{verbatim}

\usepackage{caption}
\usetikzlibrary{shapes}
\usepackage{subcaption}
\usepackage{geometry}
\usepackage{mathrsfs} 
\usepackage{array}
\usepackage{multirow}
\usepackage{booktabs} 
\usepackage{adjustbox}
\usepackage[title]{appendix}

\biboptions{sort&compress}
\hypersetup{
colorlinks=true,
linkcolor=red,
filecolor=green,
urlcolor=green,
anchorcolor=green,
citecolor=cyan,
pdftitle={Overleaf Example},
pdfpagemode=FullScreen,
}
\newtheorem{remark}{Remark}
\newtheorem{corollary}{Corollary}

 \usepackage{xcolor}

 \usepackage[figcolor=white]{todonotes}

\begin{document}
	
\begin{tikzpicture}[remember picture,overlay]
\node[anchor=north east,inner sep=20pt] at (current page.north east)
{\includegraphics[scale=0.2]{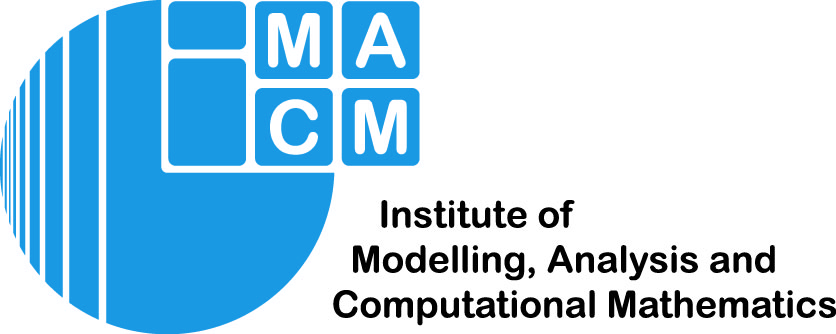}};
\end{tikzpicture}

\begin{frontmatter}

\title{A Koopman Operator Framework for Nonlinear Epidemic\\ Dynamics: Application to an SIRSD Model}

\author[BUW,MAIS]{Achraf Zinihi}
\ead{a.zinihi@edu.umi.ac.ma} 

\author[BUW]{Matthias Ehrhardt\corref{Corr}}
\cortext[Corr]{Corresponding author}
\ead{ehrhardt@uni-wuppertal.de}

\author[MAIS]{Moulay Rchid Sidi Ammi}
\ead{rachidsidiammi@yahoo.fr}

\address[BUW]{University of Wuppertal, Applied and Computational Mathematics,\\
Gaußstrasse 20, 42119 Wuppertal, Germany}

\address[MAIS]{Department of Mathematics, AMNEA Group, Faculty of Sciences and Techniques,\\
Moulay Ismail University of Meknes, Errachidia 52000, Morocco}



\begin{abstract}
We develop and analyze an SIRSD epidemic model, which extends the classical SIR framework by incorporating waning immunity and disease-induced mortality. A rigorous well-posedness analysis ensures the existence, uniqueness, positivity, and boundedness of solutions, guaranteeing the model's epidemiological feasibility. To facilitate theoretical investigations and data-driven modeling, we reformulated the system in normalized variables. 
To capture and predict complex nonlinear epidemic dynamics, we use the Koopman operator framework with extended dynamic mode decomposition (EDMD) and an epidemiologically informed dictionary of observables. 
We compare two Koopman approximations: 
one based on a minimal epidemiological dictionary and another enriched with nonlinear and cross terms. We generate synthetic data using a nonstandard finite difference (NSFD) scheme for four representative epidemics: SARS-CoV-2,
seasonal influenza, Ebola, and measles. Numerical experiments demonstrate that the Koopman-based approach effectively identifies dominant epidemic modes and accurately predicts key outbreak characteristics, including peak infection dynamics.

\end{abstract}


\journal{} 




\end{frontmatter}


\section{Introduction}\label{S1}
Mathematical epidemiology has long provided essential tools for understanding and predicting the spread of infectious diseases through compartmental models. Classical frameworks such as the (Susceptible–Infected–Recovered) SIR system, and its numerous extensions, describe population-level disease dynamics using ordinary, delay, or partial differential equations \cite{Bitsouni2024, Sekiguchi2011, Zinihi2025FDE, Lopez2025, GuerreroFlores2023, Zinihi2025S}. 
These models have been successfully applied to diverse epidemiological settings, ranging from forecasting the spread of the SARS-CoV-2 virus \cite{Puglisi2021, SotoRocha2026, Simeonov2023} to analyzing spatiotemporal outbreak patterns \cite{Zinihi2025NSFD, Chang2022, She2024}.

Recent advances have substantially broadened the scope of epidemic modeling. These novel approaches include fractional-order formulations with generalized operators and kernels \cite{Zinihi2025MM}, nonlinear SEIR extensions capturing psychological effects and discretization schemes \cite{Hoang2023}, bifurcation analyses of discrete-time SIR models \cite{Suo2025}, and actuarial frameworks linking epidemic dynamics with insurance mathematics \cite{Zinihi2025A}. Together, these studies demonstrate the versatility and cross-disciplinary impact of epidemic models. However, the nonlinearity remains a significant challenge for analytical tractability, real-time forecasting, and the design of effective intervention strategies.

The Koopman operator framework offers a promising solution to these challenges. It recasts nonlinear dynamics as linear evolutions in a higher-dimensional space of observables \cite{Brunton2022, Berry2025}. This lifting approach allows for the direct application of linear system theory for prediction, estimation, and control while remaining inherently data-driven \cite{Manzoor2023, Mezi2024}. In epidemiology, Koopman-based methods, particularly Dynamic Mode Decomposition (DMD) and its variants, have already demonstrated their ability to extract spatiotemporal patterns, forecast infection trajectories, and evaluate the effects of exogenous factors such as human mobility \cite{Proctor2015, Mustavee2022}. Extensions that incorporate control inputs, such as DMD with control (DMDc), further provide a principled framework for modeling vaccination, treatment, and behavioral interventions \cite{Proctor2018, Takeishi2017}. Recent computational advances, including tensorized and neural Koopman architectures \cite{Xu2025, Han2025}, have mitigated the curse of dimensionality. These advances enable richer observable sets and improved fidelity in large-scale epidemic systems.

These developments are particularly well-suited to extended models such as the (Susceptible-Infected-Recovered-Susceptible (again)-Deceased) SIRSD system \cite{basit2025calculating}, which must simultaneously capture nonlinear reinfection dynamics and intervention effects. Integrating SIRSD dynamics into a lifted Koopman space allows for accurate real-time forecasting, facilitates robust control under practical constraints. It also allows for  natural accommodation of regime shifts such as the emergence of new variants or policy interventions.

In this work, we develop a Koopman operator framework tailored specifically to modeling nonlinear epidemic dynamics, with a focus on the SIRSD model. Our main contributions are threefold:
First, we construct finite-dimensional Koopman approximations that faithfully capture the essential features of epidemic dynamics.
Second, we incorporate control inputs to enable the design and evaluation of potential intervention strategies.
Third, we demonstrate the predictive power and robustness of the proposed framework across diverse epidemiological scenarios.
Bridging nonlinear epidemic modeling with operator-theoretic and data-driven techniques, our framework provides a mathematically rigorous, computationally efficient, and practically relevant approach to epidemic analysis and management.

The remainder of this paper is structured as follows.
Section~\ref{S2} introduces the SIRSD epidemic model, including its formulation, well-posedness analysis, and normalization.
Section~\ref{S3} reviews the Koopman operator framework and its finite-dimensional approximation via extended dynamic mode decomposition (EDMD).
Section~\ref{S4} reformulates the normalized SIRSD model into a Koopman-ready form and analyzes its Jacobian structure.
Section~\ref{S5} presents the nonstandard finite difference (NSFD) scheme used to generate synthetic data.
Section~\ref{S6} reports numerical experiments for four representative epidemics: COVID-19, seasonal influenza, Ebola, and measles. These experiments illustrate the effectiveness of the Koopman-based approximation.
Finally, Section~\ref{S7} summarizes the main findings and discusses potential directions for future research.


\section{Proposed SIRSD Model}\label{S2}
In this section, we will introduce the mathematical framework that we use to model the transmission dynamics of infectious diseases. Specifically, we use an SIRSD compartmental model, which divides the population into four distinct health states.

\subsection{SIRSD Model Description}\label{S2.1}
We consider an extended deterministic SIR-type model that includes both \textit{loss of immunity} and \textit{disease-induced mortality}, with frequency-dependent transmission. The population is divided into four mutually exclusive compartments: susceptible $S(t)$, infected $I(t)$, recovered $R(t)$, and deceased $D(t)$. Let
\begin{equation*}
    N_L(t)=S(t)+I(t)+R(t),
\end{equation*}
denote the \textit{total living population}. 
The proposed SIRSD model is governed by the nonlinear system 
\begin{equation}\label{E2.1}
\left\{\begin{aligned}
\dot{S}(t) &= -\beta \,\frac{S(t)\,I(t)}{N_L(t)} + \omega\, R(t), \\
\dot{I}(t) &= \beta \,\frac{S(t)\,I(t)}{N_L(t)} - (\gamma + \mu)\,I(t), \\
\dot{R}(t) &= \gamma \,I(t) - \omega\, R(t),\\  
\dot{D}(t) &= \mu \,I(t),
\end{aligned}\right.
\end{equation}
supplied with the initial conditions
\begin{equation}\label{E2.2}
     S(0)=S_0>0, \quad I(0)=I_0>0, \quad R(0)=R_0\ge0, \quad D(0)=D_0\ge0,
\end{equation}
where $\beta>0$ is the transmission rate, $\gamma>0$ the recovery rate, $\mu>0$ the disease-induced mortality rate, and $\omega\ge0$ the rate of immunity loss (transition from recovered back to susceptible).
For simplicity, we neglected the birth and natural mortality rates, as they have much slower dynamics than the epidemic.

In this case, the epidemic dynamics evolve as follows: susceptible individuals become infected through effective contact with infected individuals at a rate $\beta \,\tfrac{S(t) \,I(t)}{N_L(t)}$. 
Thus, the force of infection depends on the fraction of infected individuals in the living population, causing a transition from the susceptible to the infected compartment. 
Infected individuals either recover at rate $\gamma$, moving into $R$, or die from the disease at rate $\mu$, moving into $D$. 
Recovered individuals have temporary immunity and return to the susceptible pool at a rate $\omega$; they are assumed to have permanent immunity when $\omega = 0$.
The deceased compartment accumulates disease-induced deaths and does not contribute to transmission.

The model \eqref{E2.1} is nonlinear due to the frequency-dependent bilinear term $\beta \,\tfrac{S(t)\,I(t)}{N_L(t)}$, and it also contains state-dependent feedback through the living population $N_L$ (which typically decreases as $D$ grows). These features generate more complex dynamics than the classical SIR model, prompting both analytical and numerical studies.

\subsection{Well-posedness Analysis}\label{S2.2}
This section analyzes the existence, uniqueness, positivity, and boundedness of solutions to the system~\eqref{E2.1}--\eqref{E2.2}. 
First, to verify the invariance and boundedness of the proposed model, consider the total population $N(t) = N_L(t) + D(t) = S(t) + I(t) + R(t) + D(t)$.
Differentiating and summing the system~\eqref{E2.1} yields
\begin{equation*}
     \dot{N}(t) = \dot{S}(t) + \dot{I}(t) + \dot{R}(t) + \dot{D}(t) \stackrel{\eqref{E2.1}}{=} 0.
\end{equation*}
This implies that the total population $N$ remains constant over time, i.e. $N = N_0 = S_0 + I_0 + R_0 + D_0$. Therefore, all solutions to \eqref{E2.1}–\eqref{E2.2} are bounded.
On the other hand, summing only the first three equations in \eqref{E2.1} gives
\begin{equation*}
   \dot{N}_L(t) = \dot{N}(t) - \dot{D}(t) = -\mu I(t).
\end{equation*}
When $I(t) \geq 0$, it follows that $\dot{N}_L(t) \leq 0$, meaning the living population $N_L$ decreases over time due to disease-induced mortality. Since $N_L(0) = S_0 + I_0 + R_0 > 0$, and that deaths occur continuously rather than instantaneously, $N_L(t)$ remains strictly positive for all finite $t \geq 0$.
If, hypothetically, $\dot{N}_L(t) \geq 0$ (which is mathematically possible but not realistic in this context), then $N_L$ would increase over time, but would remain strictly positive for all finite $t$ due to the positivity of the initial condition.

Let us now prove the positivity of the system solutions. 
We start with the infected class $I$. Multiply both sides of its differential equation by the negative part $I^- = \max(0, -I)$, and we get
\begin{equation*}
     \frac{1}{2} \frac{d}{dt} (I^-)^2 = \Bigl(\gamma + \mu - \beta \,\frac{S(t)}{N_L(t)}\Bigr) (I^-)^2.
\end{equation*}
Afterwords,
\begin{equation*}
   \bigl(I^-(t)\bigr)^2 = (I^-_0)^2 \exp\biggl[ \int_0^t \Bigl(\gamma + \mu - \beta \,\frac{S(y)}{N_L(y)}\Bigr) \,dy \biggr],
\end{equation*}
Since $I_0 > 0$, it follows that $I^-_0 = 0$. Therefore, $I^-(t) = 0$ for all $t \geq 0$, implying that $I(t) \geq 0$. 
Since $\dot{D}(t) = \mu I(t) \geq 0$ and the initial condition satisfies $D_0 \geq 0$, it follows that $D(t) \geq 0$ for all finite $t \geq 0$.
To demonstrate the nonnegativity of $R$ and $S$ individually, we use the following methodology outlined in \cite[Page~7]{Zinihi2025A} and consider
\begin{equation*}
     \dot{R}\big|_{R=0} = \gamma I \geq 0, \quad \text{and} \quad \dot{S}\big|_{S=0} = \omega R \geq 0.
\end{equation*}
Therefore, the vector field on the boundaries of the positive orthant $\mathbb{R}_+^2$ points inward or is tangent to the boundary. This ensures that solutions remain in $\mathbb{R}_+^2$.

We use standard results from the theory of  ordinary differential equations (ODE) to prove the existence and uniqueness of solutions to the system of equations~\eqref{E2.1}--\eqref{E2.2}. This is supported by the boundedness and positivity of the solutions.
The right-hand side of the system is continuously differentiable and thus locally Lipschitz continuous in the positive orthant (excluding division by zero since $N_L(t) > 0$ for all finite $t$).
According to the Picard–Lindelöf theorem, this guarantees the local existence and uniqueness of solutions with given nonnegative initial conditions.
Furthermore, the positivity of the state variables $S$, $I$, $R,$ and $D$, together with the invariance of a bounded region, ensures that the solutions will remain positive and bounded for all $t \geq 0$. This precludes finite-time blow-up and allows us to extend local solutions to global ones.
Therefore, the following corollary summarizes this section
\begin{corollary}\label{C1}
Given any initial data satisfying \eqref{E2.2} and positive parameters $\beta, \gamma, \mu, \omega$, the system~\eqref{E2.1} admits a unique, global, positive, and bounded solution for all finite $t \geq 0$.
\end{corollary}

\subsection{Model Normalization}\label{S2.3}
Our proposed model reveals that the variables $S$, $I$, and $R$ evolve independently of $D$. Therefore, it is sufficient to focus on the \textit{reduced system}
\begin{equation*}
\left\{
\begin{aligned}
\dot{S}(t) &= -\beta \,\frac{S(t)\,I(t)}{N_L(t)} + \omega\, R(t), \\
\dot{I}(t) &= \beta \,\frac{S(t)\,I(t)}{N_L(t)} - (\gamma + \mu)\, I(t), \\
\dot{R}(t) &= \gamma\, I(t) - \omega \,R(t),
\end{aligned}
\right.
\end{equation*}
where $D(t) = N - N_L(t)$.  

To simplify the analysis and prepare the model for the Koopman operator-theoretic framework, we introduce dimensionless variables representing proportions of the total population
\begin{equation*}
   s(t) = \frac{S(t)}{N}, \quad i(t) = \frac{I(t)}{N}, \quad r(t) = \frac{R(t)}{N}, \quad d(t) = \frac{D(t)}{N}.
\end{equation*}
The current living population is then
\begin{equation*}
   N_L(t) = S(t) + I(t) + R(t) = N - D(t) = N \big( 1 - d(t) \big),
\end{equation*}
which decreases over time due to disease-induced mortality. 
Dividing each equation by the constant $N$ yields the normalized system
\begin{equation}\label{E2.3}
\left\{
\begin{aligned}
\dot{s}(t) &= -\beta \frac{s(t) i(t)}{1 - d(t)} + \omega r(t), \\
\dot{i}(t) &= \beta \frac{s(t) i(t)}{1 - d(t)} - (\gamma + \mu) i(t), \\
\dot{r}(t) &= \gamma i(t) - \omega r(t),
\end{aligned}
\right.
\end{equation}
with initial conditions
\begin{equation}\label{E2.4}
s(0) = s_0 = \frac{S_0}{N}, \quad i(0) = i_0 = \frac{I_0}{N}, \quad r(0) = r_0 = \frac{R_0}{N}, \quad d(0) = d_0 = \frac{D_0}{N},
\end{equation}
where
\begin{equation}\label{E:d}
    d(t) = 1 - s(t) - i(t) - r(t).
\end{equation}

For the following reasons, we will focus our Koopman analysis on \eqref{E2.3}--\eqref{E2.4}:
($i$) Dimensionless variables make parameters easier to interpret and compare across different populations or datasets.
($ii$) The state variables $s(t)$, $i(t)$, $r(t)$, $d(t)$ are naturally bounded in $[0,1]$, which facilitates both theoretical analysis and numerical simulations.
($iii$) For Koopman operator analysis, normalized variables avoid scaling issues in observable functions and yield a more uniform representation of the dynamics across initial conditions.
($iv$) The factor $\frac{1}{1-d(t)}$ explicitly captures the effect of a shrinking living population on the transmission term, preserving the model’s accuracy under mortality.

\section{Koopman Operator Theory}\label{S3}
The Koopman operator framework provides a powerful alternative to traditional nonlinear analysis. It accomplishes this by shifting the dynamics from the original state space to a higher-dimensional (possibly infinite-dimensional) space of observables. 
In this new space, the evolution becomes linear.
This property allows us to use linear operator theory to predict, control, and perform spectral analysis on nonlinear systems, including epidemiological models such as the SIRSD system.

Let the state be $\mathbf{x} = (s, i, r, d)^\top$ and $F^t\colon\mathbb{R}^4 \to \mathbb{R}^4$ be the flow map associated with the system~\eqref{E2.3}, so that
\begin{equation*}
     \mathbf{x}(t) = F^t(\mathbf{x}_0), \quad t \ge 0.
\end{equation*}
For a scalar-valued observable $g\colon\mathbb{R}^4 \to \mathbb{R}$, the Koopman operator $\mathcal{K}^t$ is defined as
\begin{equation}\label{E3.1}
      (\mathcal{K}^t g)(\mathbf{x}) = g(F^t(\mathbf{x})).
\end{equation}
By definition, the Koopman operator $\mathcal{K}^t$ is linear in $g$ even though $F^t$ is generally nonlinear in the state variable $\mathbf{x}$. 
This linearity allows us to analyze nonlinear epidemic dynamics via the spectral properties of $\mathcal{K}^t$. 
However, since the Koopman operator acts on an infinite-dimensional space of observables, finite-dimensional approximations are necessary for computational purposes.

\subsection{Dictionary of Observables}\label{S3.1}
To numerically approximate $\mathcal{K}^t$, we define a finite dictionary of observables
\begin{equation*}
     \mathcal{D} = \{\psi_1, \psi_2, \dots, \psi_N\},
\end{equation*}
where each $\psi_j \colon\mathbb{R}^4 \to \mathbb{R}$ is a scalar function of the state variables.
A minimal choice for the SIRSD model~\eqref{E2.3} includes the linear observables
\begin{equation*}
\psi_1(s,i,r,d) = s, \quad \psi_2(s,i,r,d) = i, \quad \psi_3(s,i,r,d) = r, \quad \psi_4(s,i,r,d) = d.
\end{equation*}
To capture nonlinear epidemiological interactions such as the infection term $\frac{si}{1-d}$, the dictionary is enriched with polynomial and rational terms
\begin{equation*}
\psi_5(s,i,r,d) = s\, i, \quad \psi_6(s,i,r,d) = s \,r, \quad \psi_7(s,i,r,d) = i\, r, \quad \psi_8(s,i,r,d) = \frac{s \,i}{1-d}.
\end{equation*}
Depending on the desired approximation accuracy, additional higher-order terms can be included
\begin{equation*}
\psi_9(s,i,r,d) = s^2, \quad \psi_{10}(s,i,r,d) = i^2, \quad \psi_{11}(s,i,r,d) = r^2, \quad \psi_{12}(s,i,r,d) = d^2, \; \dots
\end{equation*}

In the Koopman operator framework, the choice of $\mathcal{D}$ plays a central role in capturing the nonlinear dynamics of epidemic models. 
In our SIRSD model~\eqref{E2.3}, in addition to the epidemiologically meaningful bilinear incidence term $si$, we enrich the lifted observable space with additional quadratic and cross terms such as $\psi_5$, $\psi_6$, $\psi_7$, $\psi_9$, $\psi_{10}$, $\psi_{11}$, and $\psi_{12}$. 
While most of these terms do not correspond to direct epidemiological processes, except for for instance, the term modeling new infections ($\psi_5$), they are mathematically valuable because they expand the basis functions used in the numerical approximation. 
Quadratic self-terms such as $\psi_9, \psi_{10}$ are also employed in the modeling literature, where they are used to capture nonlinear feedback or density-dependent effects (see, e.g., \cite{Verma2021}).

To assess the numerical impact of the dictionary $\mathcal{D}$, we will compare two dictionaries in our simulations: 
($\mathcal{D}_1$) a dictionary containing only the basic compartments and $\psi_8$; ($\mathcal{D}_2$) an extended dictionary including the additional quadratic and cross terms.
This comparison illustrates how adding nonlinear observables to the dictionary affects the Koopman-based approximation of epidemic dynamics. 
It is important to note that the presented terms are illustrative, and that the choice of observables in the Koopman framework is not limited to these terms.

\subsection{Data Collection}\label{S3.2}
The Koopman framework requires \textit{state trajectory data}. This data can be obtained either from:
(1) \textit{Numerical simulations (synthetic data)} of the SIRSD model~\eqref{E2.3} under given parameters $(\beta, \gamma, \mu, \omega)$ and initial conditions $(s_0, i_0, r_0, d_0)$, or
(2) \textit{Empirical epidemiological time series}, after appropriate normalization.

We discretize the trajectory at uniform sampling times $t_k = k \Delta t$, with $k = 0,1,\dots,M$, producing snapshot vectors
\begin{equation*}
    \mathbf{x}_k = \bigl(s(t_k), i(t_k), r(t_k), d(t_k)\bigr)^\top.
\end{equation*}
Each snapshot is lifted to the observable space via
\begin{equation*}
     \mathbf{y}_k = \bigl(\psi_1(\mathbf{x}_k), \psi_2(\mathbf{x}_k), \dots, \psi_N(\mathbf{x}_k)\bigr)^\top.
\end{equation*}
The resulting datasets
\begin{equation*}
   Y = (\mathbf{y}_0, \dots, \mathbf{y}_{M-1}) \in \mathbb{R}^{N \times M} \quad \text{ and } \quad Y' = (\mathbf{y}_1, \dots, \mathbf{y}_M) \in \mathbb{R}^{N \times M},
\end{equation*}
serve as the input for extended dynamic mode decomposition (EDMD) \cite{Korda2017, Jin2024}.

\subsection{Koopman Matrix Estimation}\label{S3.3}
In the finite-dimensional observable space, the Koopman dynamics are approximated by
\begin{equation*}
     \mathbf{y}_{k+1} \approx K  \mathbf{y}_k,
\end{equation*}
where $K \in \mathbb{R}^{N \times N}$ is the Koopman matrix. Thus, we approximate the action of the operator $\mathcal{K}^{\Delta t}$ by
\begin{equation*}
   Y' \approx K Y.
\end{equation*}
Moreover, the least-squares problem
\begin{equation*}
    K = \arg \min_{\tilde{K} \in \mathbb{R}^{N\times N}} \| Y' - \tilde{K} Y \|_F,
\end{equation*}
has the closed-form solution
\begin{equation}\label{E3.2}
   K = Y' Y^\star,
\end{equation}
where $Y^\star$ denotes the Moore–Penrose pseudoinverse of $Y$ and $\|\cdot\|_F$ denotes the Frobenius norm.

The spectrum of $K$ provides Koopman eigenvalues and modes, revealing growth and decay rates and oscillatory patterns in the epidemic dynamics.
Furthermore, once $K$ is computed, future state prediction is performed by iterating the linear system in the observable space.

\subsection{Finite-Dimensional Approximation via EDMD}\label{S3.4}
The EDMD is a data-driven algorithm that creates the aforementioned finite-dimensional Koopman approximation. The main steps are as follows:

\begin{itemize}
\item[$(i)$] Choose a dictionary $\mathcal{D} = \{\psi_1,\dots,\psi_N\}$ containing both linear and nonlinear epidemiologically relevant observables.

\item[$(ii)$] Collect snapshots $\{\mathbf{x}_k\}_{k=0}^M$ from simulation or real epidemic data.

\item[$(iii)$] Lift the data into the observable space to form $Y$ and $Y'$.

\item[$(iv)$] Compute $K$ using the closed-form solution \eqref{E3.2}.

\item[$(v)$] Analyze the spectral properties (i.e.\ eigenvalues and eigenvectors) to extract coherent patterns and predict dynamics.
\end{itemize}

This approach yields a linear, reduced-order surrogate of the original nonlinear SIRSD model. This surrogate is suitable for forecasting the spread of epidemics, performing stability analyses in observable spaces, and designing optimal control strategies within the framework of linear systems theory.

\section{Application to the SIRSD Model}\label{S4}
We now apply the Koopman operator and its finite-dimensional approximation via EDMD to the normalized SIRSD epidemic model~\eqref{E2.3}. 
In this section, we reformulate the system in a Koopman-ready form, discuss its Jacobian structure, and demonstrate how the EDMD methodology can be tailored to capture the epidemiological dynamics.

\subsection{Koopman-ready Form of the Normalized SIRSD Model}\label{S4.1}
Let the state be $\mathbf{x} = (s, i, r, d)^\top$ and 
\begin{equation*}
  f\bigl(\mathbf{x}(t)\bigr) =
  \begin{pmatrix}
  -\beta\,\frac{s\,i}{1-d} + \omega \,r \\
  \beta\,\frac{s\,i}{1-d} - (\gamma+\mu)\,i \\
   \gamma\, i - \omega \,r \\
    \mu\, i
\end{pmatrix}.
\end{equation*}
Then, the normalized SIRSD dynamics take the form
\begin{equation}\label{E4.1}
\dot{\mathbf{x}}(t) = f\bigl(\mathbf{x}(t)\bigr),
\end{equation}
where the initial state is given by
\begin{equation*}
   \mathbf{x}(0) = (s_0, i_0, r_0, d_0)^\top.
\end{equation*}
The admissible state space is the epidemic simplex:
\begin{equation*}
   \mathcal{S} = \bigl\{(s,i,r,d) \in \mathbb{R}^4_{+} \ : \ s+i+r+d = 1, \ d<1\bigr\},
\end{equation*}
ensuring $(1-d)^{-1}$ is well-defined for all finite $t$.

\begin{remark}
This formulation is Koopman-ready: the state variables $(s,i,r,d)$ evolve according to smooth, well-defined functions on $\mathcal{S}$, which allows for lifting into a finite-dimensional observable space.
\end{remark}

Following the EDMD procedure in Section~\ref{S3}, we examine two dictionaries of observables to evaluate the effect of dictionary design on the Koopman approximation of the SIRSD dynamics.
The first, minimal dictionary,
\begin{equation*}
   \mathcal{D}_1 = \Bigl\{s, \ i, \ r, \ d, \ \frac{s\,i}{1-d}\Bigr\},
\end{equation*}
contains only the basic compartments together with the nonlinear infection term $\frac{s\i}{1-d}$.
The second, extended dictionary,
\begin{equation*}
    \mathcal{D}_2 = \Bigl\{s, \ i, \ r, \ d, \ s\,i, \ s\,r, \ i\,r, \frac{s\,i}{1-d}, \ s^2, \ i^2, \ r^2, \ d^2 \Bigr\},
\end{equation*}
enriches  the first dictionary $\mathcal{D}_1$ by including cross terms and quadratic self-terms.
This comparison allows us to evaluate the influence of extending the dictionary with nonlinear observables on the numerical representation of epidemic dynamics.

\subsection{Jacobian Structure}\label{E4.2}
The Jacobian $\mathcal{J}$ of system \eqref{E4.1} at $\mathbf{x}$ is given by
\begin{equation*}
\mathcal{J}(\mathbf{x}) =
\begin{pmatrix}
-\beta \,\frac{i}{1-d} & -\beta \,\frac{s}{1-d} & \omega & -\beta \,\frac{s\,i}{(1-d)^2} \\
\beta \,\frac{i}{1-d} & \beta\, \frac{s}{1-d} - (\gamma+\mu) & 0 & \beta \,\frac{s \,i}{(1-d)^2} \\
0 & \gamma & -\omega & 0 \\
0 & \mu & 0 & 0
\end{pmatrix}.
\end{equation*}
Terms involving $(1-d)^{-1}$ and $(1-d)^{-2}$ naturally arise from the standard incidence scaling with a decreasing living population size. 
The resulting Jacobian plays an important role in several contexts: it is used for local linearization near equilibria, for constructing the Koopman generator as the continuous-time analogue of the Koopman operator, and for performing sensitivity analysis with respect to the parameters $\beta$, $\gamma$, $\mu$, and $\omega$.

\section{Numerical Scheme}\label{S5}
We employ a nonstandard finite difference (NSFD) approach to discretize the normalized SIRSD system~\eqref{E2.3}--\eqref{E2.4}. 
Originally proposed by Mickens~\cite{Mickens1993}, this scheme is particularly effective for biological and epidemiological models since it preserves essential qualitative features of the underlying continuous dynamics, such as positivity, boundedness, and dynamic consistency.

Standard numerical schemes, such as Runge-Kutta methods or the explicit/implicit Euler method, often fail to maintain these properties and may generate nonphysical artifacts, including negative compartment sizes or spurious equilibria \cite{Ehrhardt2013, Costa2024, Maamar2023, Zinihi2025NSFD}. 
For example, as demonstrated in \cite{Zinihi2025NSFD, Maamar2023}, the traditional finite difference method may yield numerically negative exposed populations under specific circumstances. In contrast, the NSFD scheme consistently preserves positivity and ensures qualitative realism.

In the present framework, two fundamental requirements are maintaining the positivity of the state variables and ensuring that the  compartmental proportions remain bounded within $[0,1]$ (as established in Section~\ref{S2}).
Furhermore, the NSFD construction guarantees dynamic consistency, meaning that key thresholds and long-term behaviors of the discrete model replicate those of the continuous formulation. These properties make the method especially suitable for the SIRSD dynamics, in which mortality and immunity loss jointly shape epidemic trajectories. Finally, the resulting NSFD discretization provides synthetic data that can be directly used within the proposed Koopman framework.

In the NSFD framework, the time derivative $\dot{u}(t)$ of a generic state variable $u$ is approximated using a nontrivial denominator function $\varphi(k)$, where $k=\Delta t$ is the time step
\begin{equation*}
  \dot{u}\bigg|_{t=t_n} \approx \frac{u^{n+1} - u^n}{\varphi(k)},
\end{equation*}
with $\varphi(k) > 0$ and $\varphi(k) = k + \mathcal{O}(k^2)$. One common choice that guarantees positivity preservation is $\varphi(k) = \tfrac{e^{\eta k} - 1}{\eta}$, where $\eta$ is the natural mortality rate, as discussed in \cite{Zinihi2025NSFD}.
Let $k$ be the time step size, and denote $s^n \approx s(t_n)$, $i^n \approx i(t_n)$, and $r^n \approx r(t_n)$. 
The NSFD discretization of the SIRSD model~\eqref{E2.3} becomes
\begin{equation}\label{E:NSFD}
\left\{\begin{aligned}
\frac{s^{n+1} - s^n}{\varphi(k)} &= - \beta \,\frac{s^{n+1} i^n}{1 - d^n} + \omega \,r^{n}, \\
\frac{i^{n+1} - i^n}{\varphi(k)} &= \beta \,\frac{s^{n+1} i^n}{1 - d^n} - (\gamma + \mu)\, i^{n+1}, \\
\frac{r^{n+1} - r^n}{\varphi(k)} &= \gamma \,i^{n+1} - \omega \,r^{n+1}, 
\end{aligned}\right.
\end{equation}
%
where the mortality proportion $d^n \approx d(t_n)$, defined in \eqref{E:d}, is updated as
\begin{equation*} 
    \frac{d^{n+1} - d^n}{\varphi(k)} = \mu\, i^{n+1} + \omega\,\frac{r^{n+1} - r^n}{\varphi(k)}, \quad n=0,1,2\dots.
\end{equation*}
%
This scheme follows Mickens’s rules \cite{Mickens1993};
\begin{itemize}
\item Nonlinear terms such as the quadratic infection term $\beta \tfrac{s\, i}{1-d}$ are discretized non-locally, i.e., using a mix of time levels (e.g., $s^{n+1} \,i^n$) to maintain positivity;

\item The discrete derivative uses a nonlinear denominator function $\varphi(k)$ that reflects the asymptotic behavior of the system;

\item The scheme is explicitly solvable in a sequential manner, with each variable updated in a given order.
\end{itemize}
The proposed NSFD scheme~\eqref{E:NSFD} can be rewritten and solved sequentially as follows
\begin{equation}\label{E5.1}
\left\{\begin{aligned}
   s^{n+1} &= \frac{s^n + \omega \, r^{n+1}\varphi(k)}{1 + \beta \,\frac{i^n}{1 - d^n}\,\varphi(k)}, \\
   i^{n+1} &= \frac{i^n + \beta \,\frac{s^{n+1} \, i^n}{1 - d^n}\,\varphi(k)}{1 + (\gamma + \mu)\,\varphi(k)}, \\
    r^{n+1} &= \frac{r^n + \gamma \,i^{n+1}\,\varphi(k)}{1 + \omega\,\varphi(k)}, \\
    d^{n+1} &= 1 - s^{n+1} - i^{n+1} - r^{n+1}.
\end{aligned}\right.
\end{equation}

We will briefly discuss the discretization of the nonlinear incidence term $\beta \tfrac{s i}{1-d}$. 
Specifically,  in the first two equations of \eqref{E5.1}, the infection term is discretized as $\beta\tfrac{s^{n+1}\, i^n}{1 - d^n}$, rather than $\beta\tfrac{s^n\,i^n}{1 - d^n}$ or $\beta\tfrac{s^{n+1}\,i^{n+1}}{1 - d^{n+1}}$. 
The guiding principle here is to evaluate exactly one factor at the new time level, specifically, the variable whose time derivative appears in the equation (in this case, $s$ or $i$). 
This semi-implicit treatment ensures the positivity of the numerical solution while maintaining the explicit solvability of the scheme. 
This strategy aligns with the general rules of NSFD methods and has been successfully applied in similar epidemiological contexts (see e.g., \cite{Ehrhardt2013, Mickens1993, Zinihi2025NSFD, Zinihi2025A}).

The NSFD scheme \eqref{E5.1} is designed so
that all compartments $s^n$, $i^n$, $r^n$, $d^n$ remain non-negative and the discrete total normalized population 
$s^n+i^n+r^n+d^n$ remains constant at 1, 
provided that the initial data satisfy this normalization condition. 
In what follows, the time step $k$ used in the NSFD scheme~\eqref{E5.1} is set to $k = 0.1$, and the denominator function is chosen as $\varphi(k) = \tfrac{e^{\eta k} - 1}{\eta}$, with $\eta\to0$ (i.e., no natural mortality), simplifying to $\varphi(k) = k$.

\section{Numerical Experiments}\label{S6}
For the simulations, we is discretized the SIRSD model~\eqref{E2.3}--\eqref{E2.4} using the NSFD scheme~\eqref{E5.1} to generate synthetic data suitable for Koopman operator-based approximation. 
The model parameters are chosen based on published epidemiological data for four representative epidemics: COVID-19, seasonal influenza, Ebola, and measles \cite{JHU2020, CDCCOVID, CDCFlu, ECDCFlu, CDCEbola2014, WHOEbola2014, CDCMeasles, CDCPinkBook}.
These parameter values are summarized in Table~\ref{Tab1} 
and are used to analyze the system's dynamics under four distinct scenarios.
This approach allows for a systematic evaluation of how the model responds to variations in transmission, recovery, and mortality rates and demonstrates the effectiveness of the proposed Koopman framework across various synthetic data cases.

\begin{table}[H]
\centering
\setlength{\tabcolsep}{0.3cm}
\caption{Epidemiological parameters for the SIRSD model applied to selected epidemics}\label{Tab1}
\adjustbox{max width=\textwidth}{
\begin{tabular}{cccccc}
\hline
\textbf{Epidemic} & $\beta$ & $\gamma$ & $\mu$ & $\omega$ & \textbf{References} \\
\hline\hline
COVID-19 & \multirow{2}{*}{0.3 – 0.6} & \multirow{2}{*}{0.07 – 0.1} & \multirow{2}{*}{0.005 – 0.01} & \multirow{2}{*}{0 – 0.02} & \multirow{2}{*}{\cite{JHU2020, CDCCOVID}} \\
(Wuhan 2020) & & & & &\\
\hline
Seasonal Influenza & 0.4 – 0.8 & 0.2 – 0.33 & 0.0001 – 0.001 & 0.1 – 0.3 & \cite{CDCFlu, ECDCFlu} \\
\hline
Ebola & \multirow{2}{*}{0.18 – 0.25} & \multirow{2}{*}{0.086 – 0.1} & \multirow{2}{*}{0.35 – 0.5} & \multirow{2}{*}{0} & \multirow{2}{*}{\cite{CDCEbola2014, WHOEbola2014}} \\
(West Africa 2014) & & & & &\\
\hline
Measles & \multirow{2}{*}{1.5 – 3.5} & \multirow{2}{*}{0.1 – 0.14} & \multirow{2}{*}{0.0001 – 0.002} & \multirow{2}{*}{0} & \multirow{2}{*}{\cite{CDCMeasles, CDCPinkBook}}\\
(Pre-vaccine era) & & & & &\\
\hline
\end{tabular}
}
\end{table}

We initialize the system with normalized compartments (proportions) that sum to one:
\begin{equation*}
    s(0)+i(0)+r(0)+d(0)=1, \ \text{ with } \ s(0)=1-i(0),\ r(0)=0,\ d(0)=0,
\end{equation*}
with a small infectious seed $i(0)=10^{-1}$ (one per $10$ individuals), representing the early stage of the outbreak.

\begin{figure}[H]
\centering
\includegraphics[width=\textwidth]{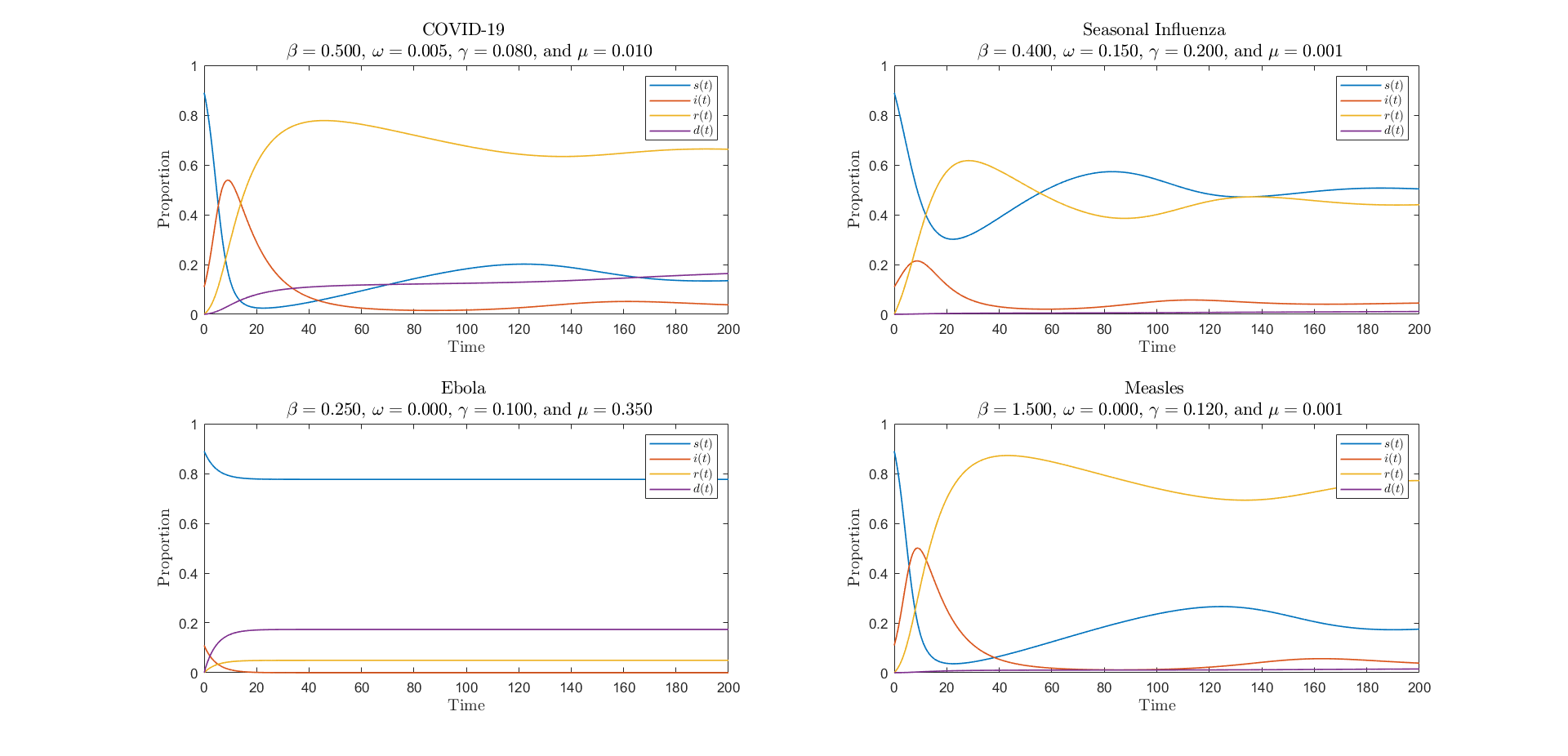}\vspace{-0.5cm}
\caption{Dynamics of the SIRSD epidemic model~\eqref{E2.3}--\eqref{E2.4} generated via the NSFD scheme~\eqref{E5.1}.}\label{F1}
\end{figure}

Figure~\ref{F1} presents numerical simulations of the SIRSD model~\eqref{E2.3}--\eqref{E2.4} for the four representative infectious diseases.  Each subplot shows how the proportion of susceptible $s(t)$, infected $i(t)$, recovered $r(t)$, and deceased $d(t)$ individuals in the population changes over time.\\
For example, for the case of the top-left subplot representing the spread of the SARS-CoV-2 virus (COVID-19)
with parameters set to $\beta = 0.5$, $\omega = 0.005$, $\gamma = 0.08$, and $\mu = 0.01$, infections rise rapidly before declining as recoveries peak, while mortality gradually increases. The susceptible population decreases markedly before reaching a stable level.\\
For seasonal influenza (top-right), with $\beta = 0.4$, $\omega = 0.15$, $\gamma = 0.2$, and $\mu = 0.001$, infections exhibit a moderate peak and faster recovery, leading to oscillatory dynamics between susceptible and recovered individuals due to reinfection.\\
For measles (bottom-right), with $\beta = 1.5$, $\omega = 0$, $\gamma = 0.12$, and $\mu = 0.001$, the model predicts a sharp outbreak with a high infection peak, followed by substantial recovery. At the same time, the proportion of susceptible individuals declines steeply and stabilizes.\\
In contrast, Ebola (bottom-left), with $\beta = 0.25$, $\omega = 0$, $\gamma = 0.1$, and $\mu = 0.35$, the model predicts relatively low infection levels, but a sharp increase in deaths, while most of the population remains susceptible.

These simulations demonstrate how disease-specific transmission rates, waning immunity, recovery rates, and mortality rates influence epidemic dynamics. They also generate synthetic ground-truth data for the subsequent Koopman operator learning framework.

\begin{figure}[H]
\centering
\begin{subfigure}[b]{\textwidth}
\centering
\includegraphics[width=\textwidth]{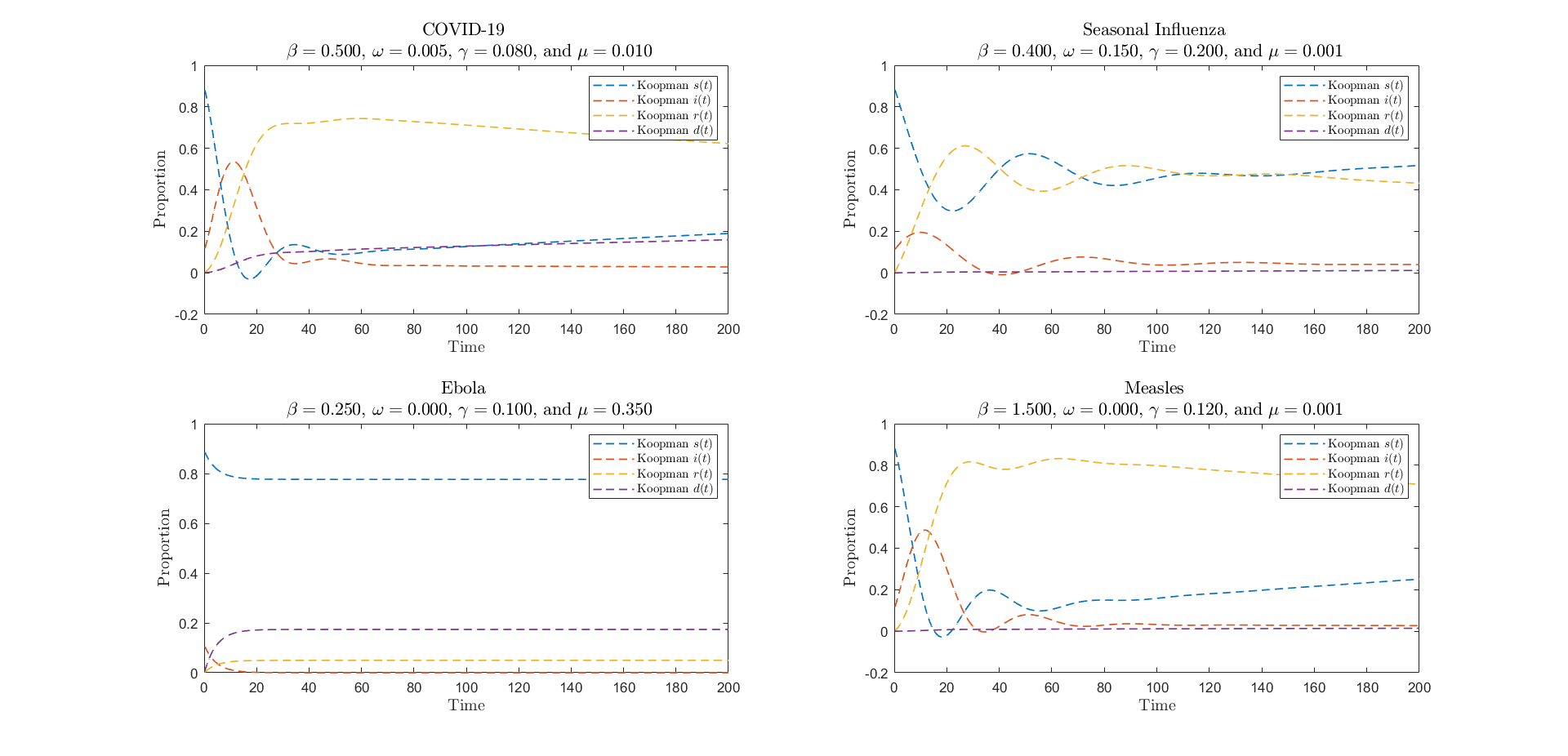}\vspace{-0.5cm}
\caption{Koopman validation with minimal dictionary $\mathcal{D}_1$.}\label{F2.1}
\end{subfigure}
\begin{subfigure}[b]{\textwidth}
\centering
\includegraphics[width=\textwidth]{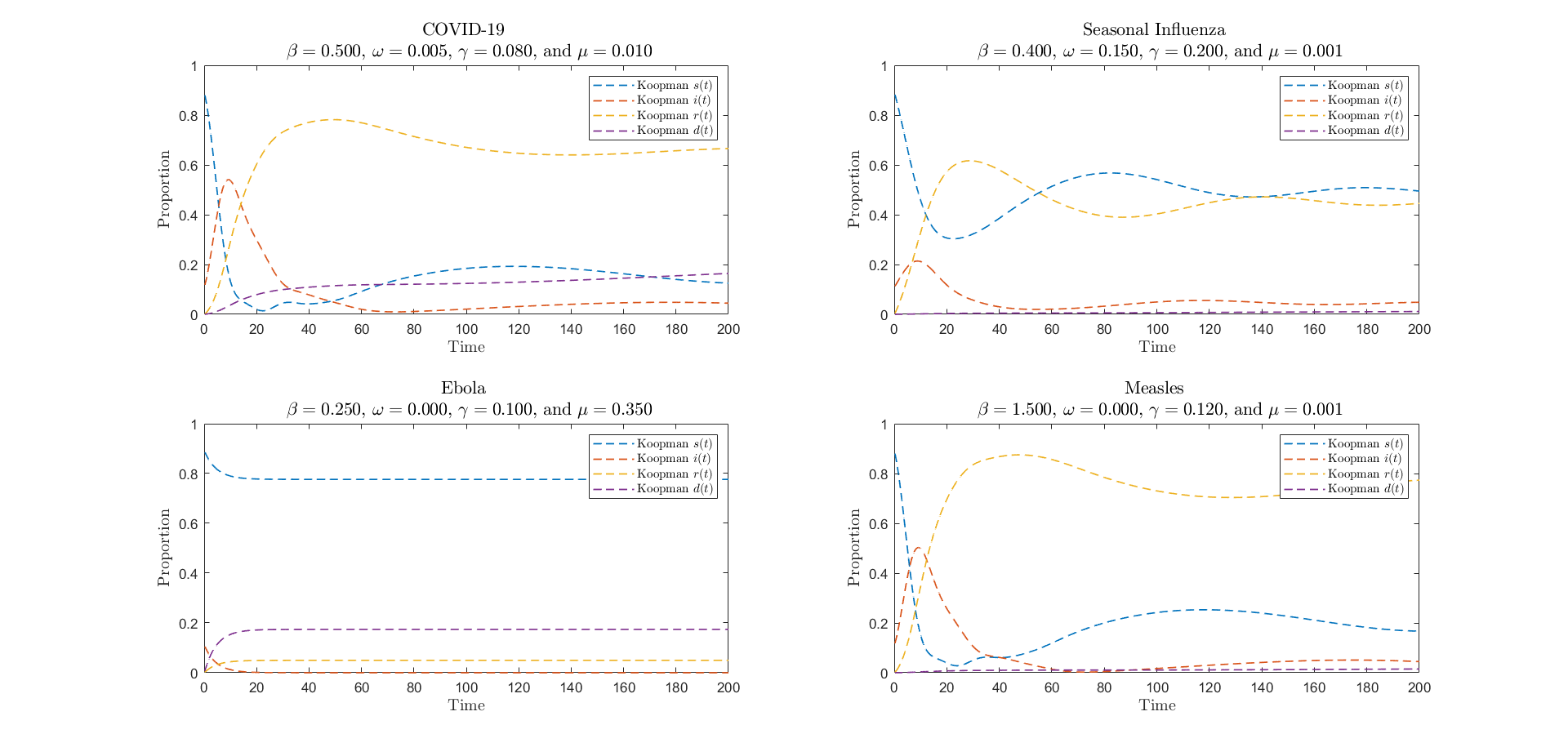}\vspace{-0.5cm}
\caption{Koopman validation with enriched dictionary $\mathcal{D}_2$.}\label{F2.2}
\end{subfigure}
\caption{Koopman-based reconstruction of SIRSD dynamics from synthetic data.}\label{F2}
\end{figure}

Figure~\ref{F2} presents the Koopman-based reconstructions of the SIRSD epidemic dynamics using EDMD with two different dictionaries of observables, $\mathcal{D}_1$ and $\mathcal{D}_2$, as mentioned in Section~\ref{S4}. 
The same epidemiological parameters as in Figure~\ref{F1} are used for the four case studies.

Subfigure~\ref{F2.1}: Koopman validation with the minimal dictionary $\mathcal{D}_1$. 
While the main epidemic trends are reproduced, the approximations display certain deficiencies, such as negative values in some compartments (notably $s(t)$ for COVID-19, $i(t)$ for seasonal influenza, and both $s(t)$ and $i(t)$ for measles). These artifacts highlight the limited representational power of the minimal dictionary.

Subfigure~\ref{F2.2}: Koopman validation with the enriched dictionary $\mathcal{D}_2$. The reconstructions more closely match the synthetic trajectories of Figure~\ref{F1}, successfully capturing infection peaks, oscillatory dynamics, and mortality accumulation without introducing unphysical negative solutions.

\begin{figure}[H]
\centering
\begin{subfigure}[b]{\textwidth}
\centering
\includegraphics[width=\textwidth]{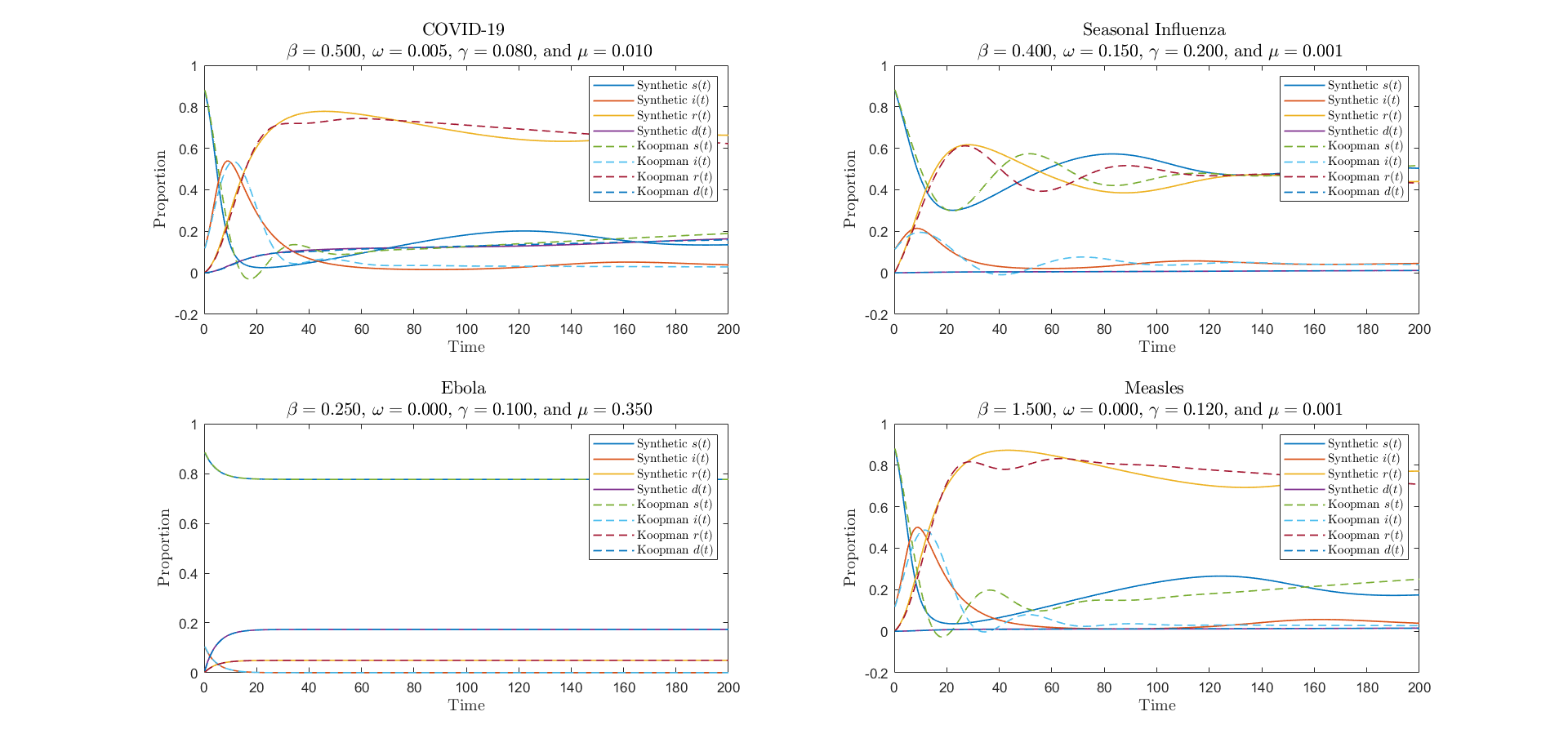}\vspace{-0.5cm}
\caption{Koopman validation with minimal dictionary $\mathcal{D}_1$.}\label{F3.1}
\end{subfigure}
\begin{subfigure}[b]{\textwidth}
\centering
\includegraphics[width=\textwidth]{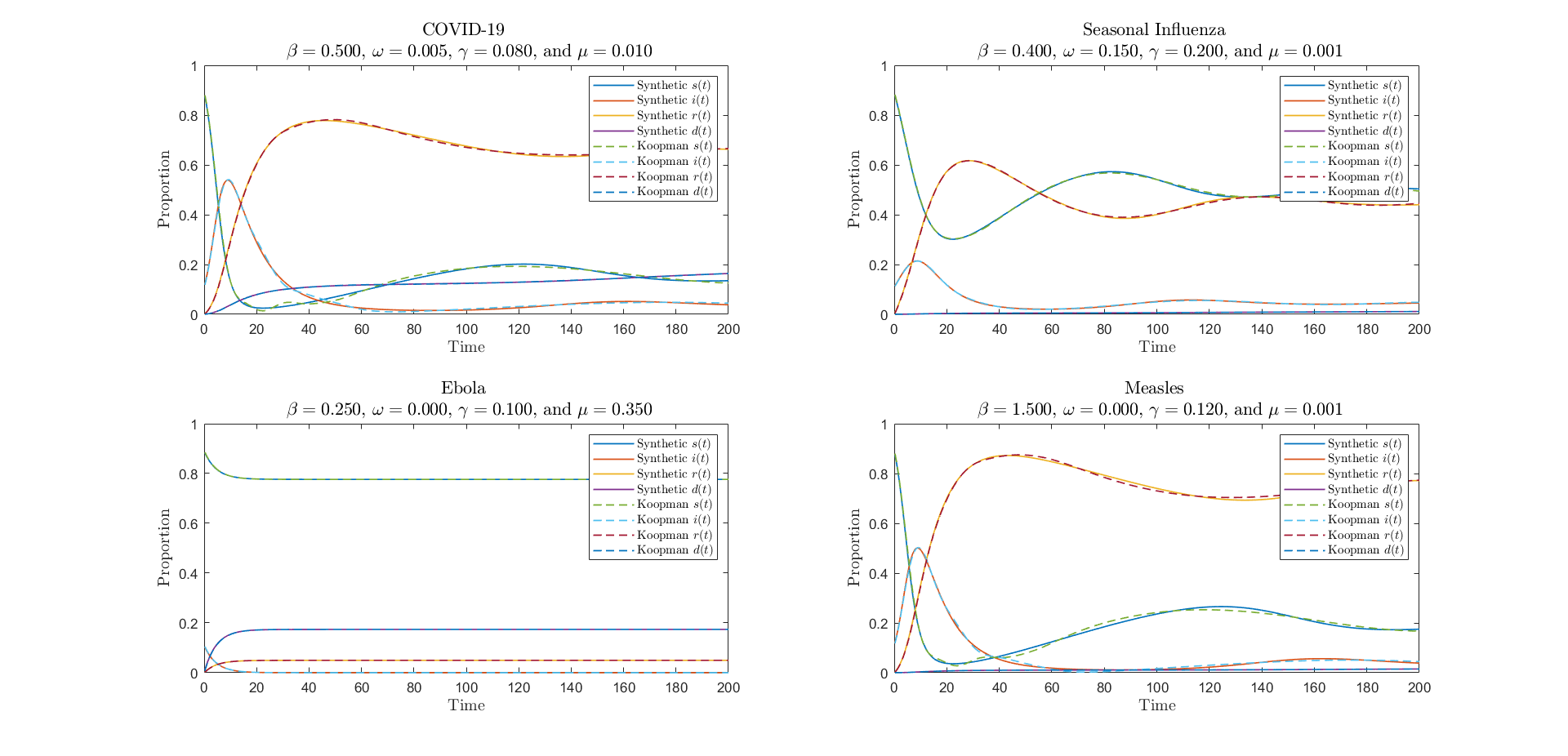}\vspace{-0.5cm}
\caption{Koopman validation with enriched dictionary $\mathcal{D}_2$.}\label{F3.2}
\end{subfigure}
\caption{Overlay comparison between synthetic SIRSD trajectories and Koopman reconstructions.}\label{F3}
\end{figure}

Overall, the enriched dictionary $\mathcal{D}_2$ demonstrates superior predictive capability compared to the minimal dictionary $\mathcal{D}_1$. 
Figure~\ref{F3} provides a direct comparison by overlaying the original synthetic data (Figure~\ref{F1}) with the two Koopman reconstructions, which further illustrate this improvement.

Figure~\ref{F3} combines the synthetic SIRSD simulations from Figure~\ref{F1} with the Koopman-based reconstructions from Figure~\ref{F2}, enabling a direct visual comparison between the ground-truth NSFD trajectories (solid lines) and the Koopman approximations (dashed lines).

Subfigure~\ref{F3.1} shows the results obtained using the minimal dictionary $\mathcal{D}_1$. 
While the main epidemic patterns are recovered, significant discrepancies appear, including mismatched oscillations and unphysical negative values in certain compartments that are physically impossible. 
These discrepancies highlight the limited representational capacity of $\mathcal{D}_1$.

\begin{figure}[H]
\centering
\includegraphics[width=\textwidth]{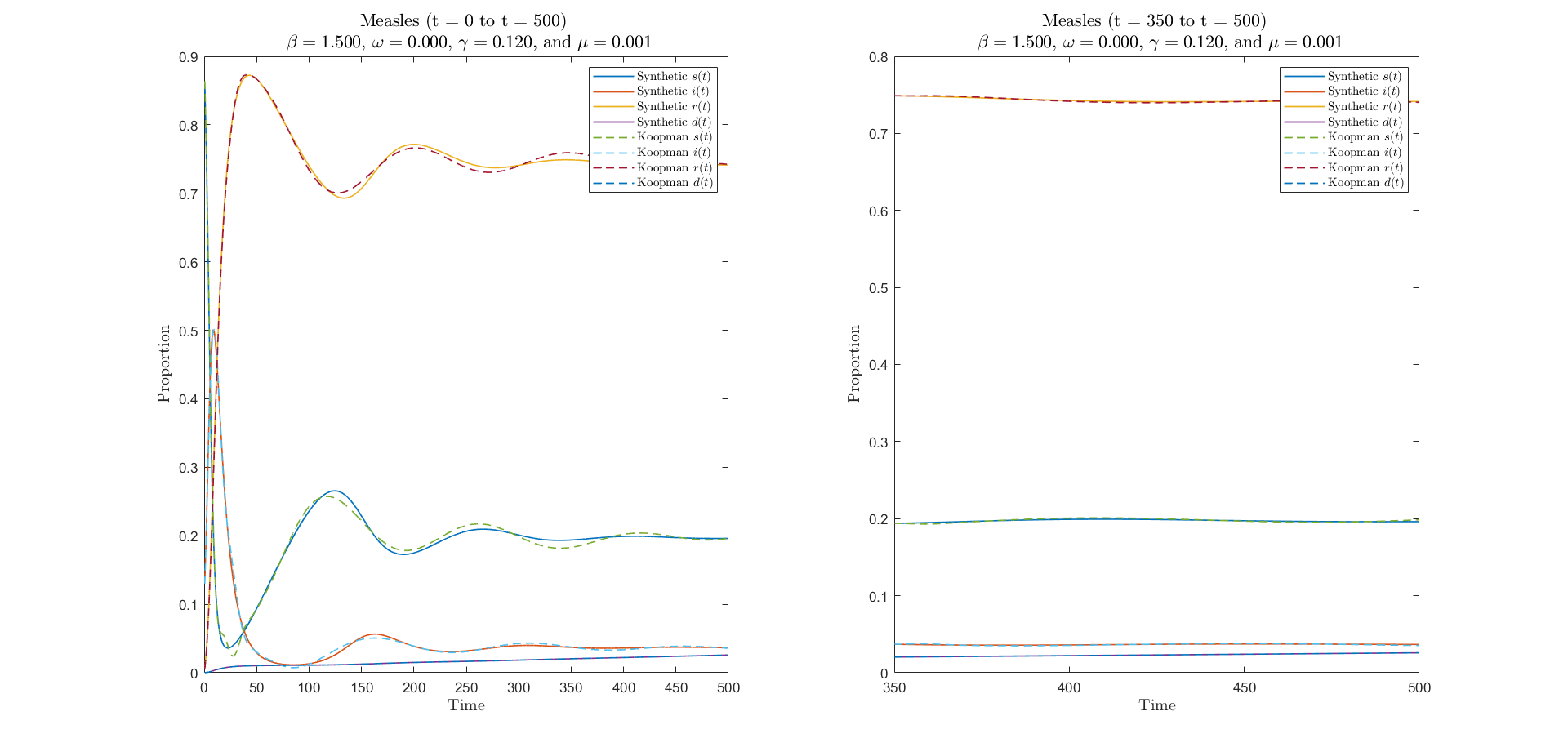}\vspace{-0.5cm}
\caption{Extended-time validation of Koopman convergence for the Measles outbreak using 
$\mathcal{D}_2$.}\label{F4}
\end{figure}

Subfigure~\ref{F3.2} depicts the enriched dictionary $\mathcal{D}_2$ and shows that the main epidemic patterns are recovered with fewer discrepancies. The Koopman reconstructions exhibit excellent agreement with the synthetic trajectories across the three diseases COVID-19, influenza, and Ebola, successfully reproducing infection peaks, oscillatory dynamics, and mortality accumulation. 
For measles, however, a small deviation emerges near the final simulation time $t=200$ in the susceptible and recovered populations, indicating slower convergence of the Koopman approximation in this highly infectious regime. 
This deviation motivates further investigation over extended time horizons, as shown in Figure~\ref{F4}.

Figure~\ref{F3} demonstrates that the enriched dictionary $\mathcal{D}_2$ provides a more robust and accurate Koopman approximation of the nonlinear SIRSD dynamics than the minimal dictionary $\mathcal{D}_1$.

Figure~\ref{F4} examines the long-term Koopman approximation of the measles case using the enriched dictionary $\mathcal{D}_2$. 
The left panel shows the entire time window $t \in [0,500]$, while the right panel focuses on the interval $t \in [350,500]$ to illustrate late-time behavior. 

Overall, the Koopman trajectories closely align with the synthetic NSFD data, successfully capturing the sharp initial outbreak, subsequent oscillations, and long-term stabilization of the system. 
The zoomed-in view confirms that the small discrepancies observed near $t=200$ in Figure~\ref{F3} dissipate over longer time horizons, with the Koopman approximation converging toward the synthetic epidemic dynamics. 

This indicates robust the Koopman approach's robust long-term forecasting capability when sufficient temporal data are available for training. 
Although the enriched dictionary $\mathcal{D}_2$ includes bilinear, quadratic, and cross terms that may not correspond directly to epidemiological mechanisms, it provides mathematically and numerically meaningful observables that enhance the representational power of the Koopman operator. This ensures improved predictive accuracy.

\section{Conclusion and Perspectives}\label{S7}
In this work, we introduced and examined an SIRSD epidemic model that considers reinfection and disease-induced mortality. 
This makes the model suitable for pathogens where partial immunity and fatality rates play a central role.
We proved the model's well-posedness, reformulated it in normalized variables, and developed an NSFD scheme to generate high-quality synthetic data.
Based on these findings, we applied the Koopman operator framework with EDMD to obtain a finite-dimensional linear approximation of the epidemic dynamics.
To evaluate the impact of dictionary design, we compared two sets of observables: a minimal epidemiological set $\mathcal{D}_1$, and an extended set $\mathcal{D}_2$ enriched with nonlinear and cross terms.
Numerical experiments conducted on four case studies (Covid-19, seasonal influenza, Ebola, and measles) show that, while the minimal dictionary $\mathcal{D}_1$ captures the essential dynamics, the enriched dictionary $\mathcal{D}_2$ yields more accurate Koopman simulations of the SIRSD system. This is particularly evident in its ability to reproduce nonlinear effects and forecast epidemic peaks. These results confirm the Koopman-SIRSD approach's capacity to capture dominant modes and provide operator-theoretic insights into disease spread.

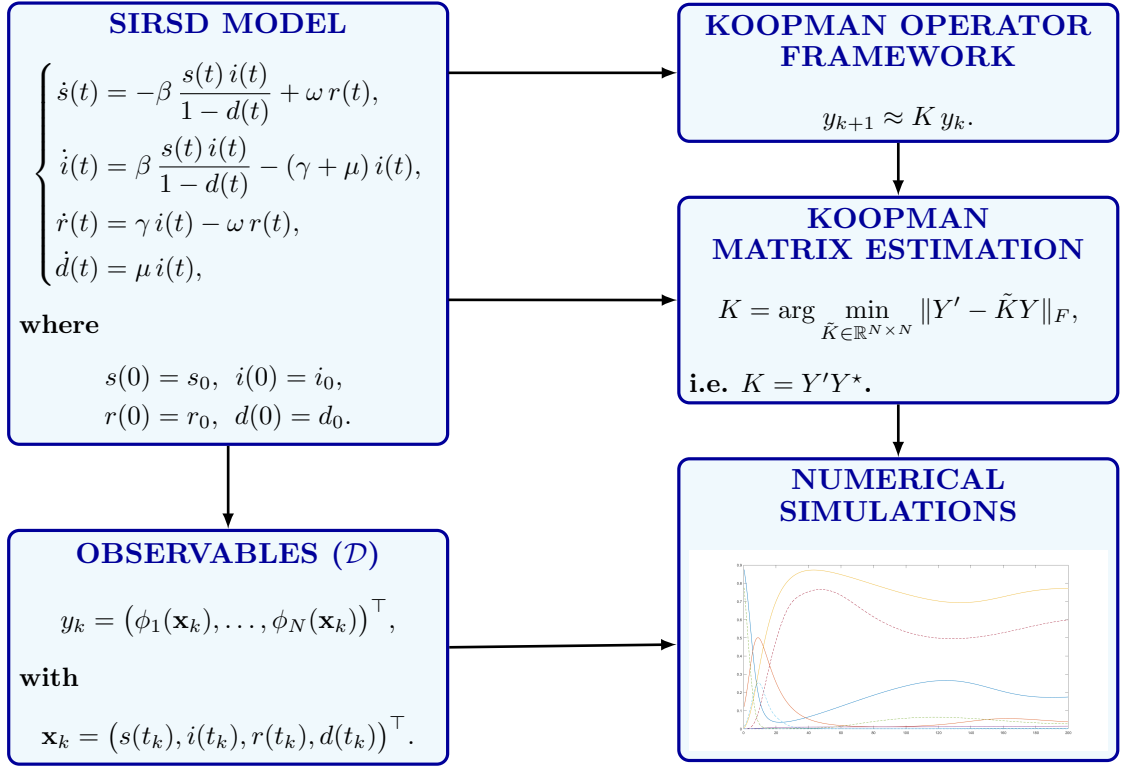
\begin{figure}[H]
\centering
\adjustbox{max width=\textwidth}{
\begin{tikzpicture}[node distance=3cm]
\node (A) [draw=blue!60!black, very thick, rounded corners,
font=\bfseries, text width=5.5cm, fill=cyan!5] 
{\small\centering{\color{blue!60!black}SIRSD MODEL} \footnotesize
$$
\left\{\begin{aligned}
\dot{s}(t) &= -\beta\,\frac{s(t)\, i(t)}{1-d(t)} + \omega\, r(t), \\
\dot{i}(t) &= \beta\,\frac{s(t)\,i(t)}{1-d(t)} - (\gamma+\mu)\,i(t), \\
\dot{r}(t) &= \gamma \,i(t) - \omega\,r(t), \\
\dot{d}(t) &= \mu \,i(t),
\end{aligned}\right.
$$
\raggedright {where}
$$
\begin{aligned}
s(0) &= s_0, \ i(0) = i_0, \\ 
r(0) &= r_0, \ d(0) = d_0.
\end{aligned}
$$
};
\node (B) [draw=blue!60!black, very thick, rounded corners,
font=\bfseries, text width=5.5cm, fill=cyan!5, below = of A, yshift = 1.9cm] 
{\small\centering{\color{blue!60!black}OBSERVABLES ($\mathcal{D}$)} \footnotesize
$$
y_k = \bigl(\phi_1(\mathbf{x}_k),\ldots,\phi_N(\mathbf{x}_k)\bigr)^\top,
$$
\raggedright {with}
$$
\mathbf{x}_k = \bigl(s(t_k), i(t_k), r(t_k), d(t_k)\bigr)^\top.
$$
};
\node (C) [draw=blue!60!black, very thick, rounded corners,
font=\bfseries, text width = 5.5cm, fill=cyan!5, right = of A, yshift = 2cm] 
{\small\centering{\color{blue!60!black}KOOPMAN OPERATOR FRAMEWORK} \footnotesize
$$
y_{k+1} \approx K\, y_k.
$$
};
\node (D) [draw=blue!60!black, very thick, rounded corners,
font=\bfseries, text width = 5.5cm, fill=cyan!5, below = of C, yshift = 2.3cm] 
{\small\centering{\color{blue!60!black}KOOPMAN\\ MATRIX ESTIMATION} \footnotesize
$$
K = \arg \min_{\tilde{K} \in \mathbb{R}^{N\times N}} \| Y' - \tilde{K} Y \|_F,
$$
\raggedright{i.e.\ $K = Y' Y^\star$.}
};
\node (E) [draw=blue!60!black, very thick, rounded corners,
font=\bfseries, text width = 5.5cm, fill=cyan!5, below = of D, yshift = 2.3cm] 
{\small\centering{\color{blue!60!black}NUMERICAL SIMULATIONS} \footnotesize
$$
\includegraphics[width=1\textwidth]{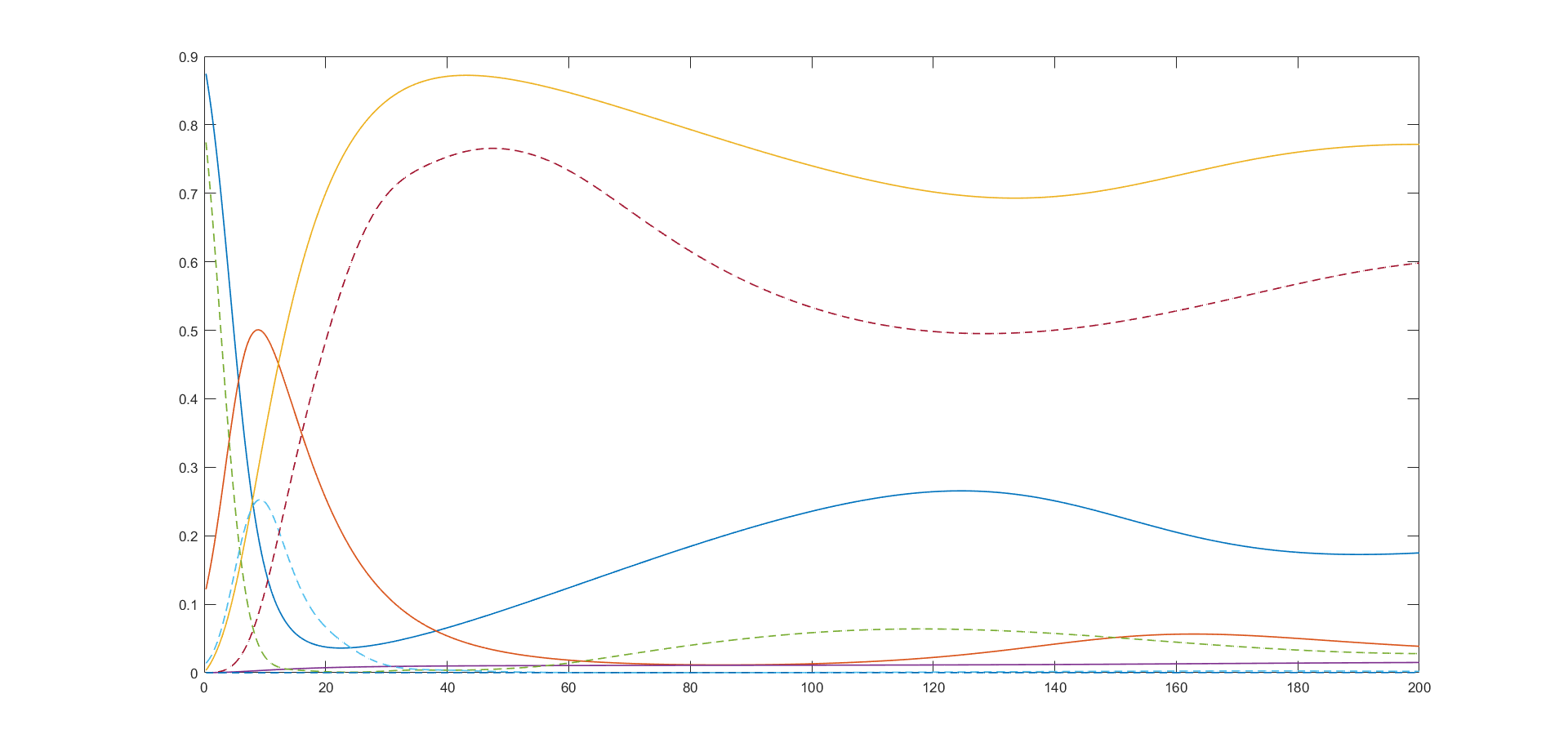}
$$
};
\draw [->, thick, >=latex, line width=1pt] (A) -- (B);
\draw [->, thick, >=latex, line width=1pt] (C) -- (D);
\draw [->, thick, >=latex, line width=1pt] (D) -- (E);
\draw [->, thick, >=latex, line width=1pt] (B.east) -- (5.9,-5.58);
\draw [->, thick, >=latex, line width=1pt] (2.9,2) -- (C.west);
\draw [->, thick, >=latex, line width=1pt] (2.9,-1.00999) -- (D.west);
\end{tikzpicture}
}
\captionof{figure}{Workflow of the SIRSD epidemic model within the Koopman operator framework.}\label{F5}
\end{figure}

In summary, Figure~\ref{F5} illustrates the overall workflow of the proposed Koopman–SIRSD framework.
Starting from the SIRSD model, we define observables and approximate the Koopman operator via EDMD. 
The resulting finite-dimensional representation enables spectral analysis and numerical simulations across different epidemic scenarios.

Beyond its immediate findings, this study reveals several promising perspectives beyond its immediate findings. The Koopman-based analysis can be extended in multiple ways:  
\begin{itemize}
\item \textbf{Control inputs:} Vaccination or treatment controls $u(t)$ can be incorporated using DMD with control (DMDc), leading to models of the form
\begin{equation*}
  Y' \approx A Y + B U,
\end{equation*}
which naturally support intervention design.  

\item \textbf{Augmented models:} Compartmental extensions such as SEIR or SIQR dynamics can be treated within the same Koopman-EDMD framework.  

\item \textbf{Spatial models:} Once discretized into ODE systems, PDE-based epidemic formulations can be lifted into Koopman space, enabling the analysis of spatiotemporal dynamics and traveling waves.

\item \textbf{Adaptive models:} Local Koopman operators may be employed to capture sudden regime changes or non-stationary behaviors in epidemic data \cite{Mezi2024}.  
\end{itemize}
This methodology is also well-suited for data-driven forecasting when applied to real epidemiological time series subject to noise and underreporting.
Finally, future research could investigate hybrid Koopman-machine learning approaches, that combine operator-theoretic structure with deep learning for robust epidemic prediction and adaptive control.

In addition to Koopman-based methodologies, there are several other research directions that can be explored. 
One such avenue is developing stochastic epidemic models, which are well-suited to capturing random effects and uncertainties arising from demographic variability or incomplete data.
Another promising direction is fractional-order epidemic models, which use derivatives of non-integer order to capture memory effects and hereditary properties in disease dynamics.
Furthermore, PDE models with linear, nonlinear, standard, and/or fractional diffusion operators can describe spatial heterogeneity and mobility-driven dynamics. 
Delayed epidemic models that incorporate incubation or reporting lags are a natural extension because time delays are known to induce richer dynamical behaviors, such as oscillations and stability switches.
Modern physics-informed neural networks (PINNs) are flexible frameworks that can be applied to stochastic, fractional, partial differential equation (PDE)-based, and delayed models. PINNs offer a unified, data-driven approach to simulating and inferring complex epidemic dynamics.

\section*{Declarations}

\subsection*{Data availability} 
All information analyzed or generated, which would support the results of this work are available in this article.
No data was used for the research described in the article.



\subsection*{Conflict of interest} 
The authors declare that there are no problems or conflicts 
of interest between them that may affect the study in this paper.





\bibliographystyle{unsrt}
\bibliography{References}

\begin{thebibliography}{10}

\bibitem{Bitsouni2024}
V.~Bitsouni, N.~Gialelis, and V.~Tsilidis.
\newblock An age‐structured {SVEAIR} epidemiological model.
\newblock {\em Mathematical Methods in the Applied Sciences},
  47(16):12460–12486, 2024.

\bibitem{Sekiguchi2011}
M.~Sekiguchi and E.~Ishiwata.
\newblock Dynamics of a discretized {SIR} epidemic model with pulse vaccination
  and time delay.
\newblock {\em Journal of Computational and Applied Mathematics},
  236(6):997–1008, 2011.

\bibitem{Zinihi2025FDE}
A.~Zinihi, M.~R. Sidi~Ammi, and D.~F.~M. Torres.
\newblock Fractional differential equations of a reaction-diffusion {SIR} model
  involving the {C}aputo-fractional time-derivative and a nonlinear diffusion
  operator.
\newblock {\em Evolution Equations and Control Theory}, 14(5):944–967, 2025.

\bibitem{Lopez2025}
S.~Lopez and N.~L. Komarova.
\newblock An optimal network that promotes the spread of an advantageous
  variant in an {SIR} epidemic.
\newblock {\em Journal of Theoretical Biology}, 605:112095, 2025.

\bibitem{GuerreroFlores2023}
S.~Guerrero‐Flores, O.~Osuna, and C.~Vargas‐De‐León.
\newblock Periodic solutions of seasonal epidemiological models with
  information‐dependent vaccination.
\newblock {\em Mathematical Methods in the Applied Sciences},
  47(4):1961–1972, 2023.

\bibitem{Zinihi2025S}
A.~Zinihi, M.~Ehrhardt, and M.~R. Sidi~Ammi.
\newblock Spatiotemporal {SEIQR} epidemic mdeling with optimal control for
  vaccination, treatment, and social measures, 2025.
\newblock arXiv preprint 2507.09328.

\bibitem{Puglisi2021}
L.~B. Puglisi, G.~S.P. Malloy, T.~D. Harvey, M.~L. Brandeau, and E.~A. Wang.
\newblock Estimation of {COVID-19} basic reproduction ratio in a large urban
  jail in the {United} {States}.
\newblock {\em Annals of Epidemiology}, 53:103–105, 2021.

\bibitem{SotoRocha2026}
M.~V.~I. Soto-Rocha, O.~Walle-García, F.~Saldaña-Jiménez,
  F.~Hernández-Cabrera, and F.J. Almaguer-Martínez.
\newblock {COVID-19} infection and recovery rates in {Mexico}: {An ARMA} time
  series model.
\newblock {\em Journal of Computational and Applied Mathematics}, 472:116781,
  2026.

\bibitem{Simeonov2023}
O.~Simeonov and C.~D. Eaton.
\newblock Modeling the drivers of oscillations in {COVID-19} data on college
  campuses.
\newblock {\em Annals of Epidemiology}, 82:40–44, 2023.

\bibitem{Zinihi2025NSFD}
A.~Zinihi, M.~Ehrhardt, and M.~R. Sidi~Ammi.
\newblock A nonstandard finite difference scheme for an {SEIQR} epidemiological
  {PDE} model, 2025.
\newblock arXiv preprint 2508.02928.

\bibitem{Chang2022}
L.~Chang, S.~Gao, and Z.~Wang.
\newblock Optimal control of pattern formations for an {SIR}
  reaction–diffusion epidemic model.
\newblock {\em Journal of Theoretical Biology}, 536:111003, 2022.

\bibitem{She2024}
G.~She and F.~Yi.
\newblock Stability and bifurcation analysis of a reaction–diffusion {SIRS}
  epidemic model with the general saturated incidence rate.
\newblock {\em Journal of Nonlinear Science}, 34(6), 2024.

\bibitem{Zinihi2025MM}
A.~Zinihi, M.~R. Sidi~Ammi, and M.~Ehrhardt.
\newblock Mathematical modeling and {Hyers-Ulam} stability for a nonlinear
  epidemiological model with {$\Phi_p$} operator and {Mittag-Leffler} kernel.
\newblock {\em Advances in Applied Mathematics and Mechanics},
  17(5):1481–1508, 2025.

\bibitem{Hoang2023}
M.~T. Hoang and M.~Ehrhardt.
\newblock A dynamically consistent nonstandard finite difference scheme for a
  generalized {SEIR} epidemic model.
\newblock {\em Journal of Difference Equations and Applications},
  30(4):409–434, 2023.

\bibitem{Suo2025}
J.~Suo and Y.~Ge.
\newblock Resonance and bifurcations of the discrete {SIR} model with unfixed
  incidence rate.
\newblock {\em Journal of Difference Equations and Applications}, page 1–27,
  2025.

\bibitem{Zinihi2025A}
A.~Zinihi, M.~Ehrhardt, and M.~R. Sidi~Ammi.
\newblock Actuarial analysis of an infectious disease insurance based on an
  {SEIARD} epidemiological model, 2025.
\newblock arXiv preprint 2508.06580.

\bibitem{Brunton2022}
S.~L. Brunton, M.~Budišić, E.~Kaiser, and J.~N. Kutz.
\newblock Modern {Koopman} theory for dynamical systems.
\newblock {\em SIAM Review}, 64(2):229–340, 2022.

\bibitem{Berry2025}
T.~Berry and S.~Das.
\newblock Limits of learning dynamical systems.
\newblock {\em SIAM Review}, 67(1):107–137, 2025.

\bibitem{Manzoor2023}
W.~A. Manzoor, S.~Rawashdeh, and A.~Mohammadi.
\newblock Vehicular applications of {Koopman} operator theory -- {A} survey.
\newblock {\em IEEE Access}, 11:25917–25931, 2023.

\bibitem{Mezi2024}
I.~Mezić, Z.~Drmač, N.~Črnjarić, S.~Maćešić, M.~Fonoberova, R.~Mohr,
  A.~M. Avila, I.~Manojlović, and A.~Andrejčuk.
\newblock A {K}oopman operator-based prediction algorithm and its application
  to {COVID-19} pandemic and influenza cases.
\newblock {\em Scientific Reports}, 14(1), 2024.

\bibitem{Proctor2015}
J.~L. Proctor and P.~A. Eckhoff.
\newblock Discovering dynamic patterns from infectious disease data using
  dynamic mode decomposition.
\newblock {\em International Health}, 7(2):139–145, 2015.

\bibitem{Mustavee2022}
S.~Mustavee, S.~Agarwal, C.~Enyioha, and S.~Das.
\newblock A linear dynamical perspective on epidemiology: {I}nterplay between
  early {COVID-19} outbreak and human mobility.
\newblock {\em Nonlinear Dynamics}, 109(2):1233–1252, 2022.

\bibitem{Proctor2018}
J.~L. Proctor, S.~L. Brunton, and J.~N. Kutz.
\newblock Generalizing {K}oopman theory to allow for inputs and control.
\newblock {\em SIAM Journal on Applied Dynamical Systems}, 17(1):909–930,
  2018.

\bibitem{Takeishi2017}
N.~Takeishi, Y.~Kawahara, and T.~Yairi.
\newblock Subspace dynamic mode decomposition for stochastic {Koopman}
  analysis.
\newblock {\em Physical Review E}, 96(3), 2017.

\bibitem{Xu2025}
Z.~Xu, R.~Dai, and K.~C. Howell.
\newblock Tensor-based {Koopman} operator and its application to optimal
  control problems.
\newblock {\em Journal of Guidance, Control, and Dynamics}, page 1–13, 2025.

\bibitem{Han2025}
M.~Han and X.~Yin.
\newblock Deep neural {Koopman} operator-based economic model predictive
  control of shipboard carbon capture system.
\newblock {\em IEEE Transactions on Control Systems Technology}, page 1–16,
  2025.

\bibitem{basit2025calculating}
A.~Basit, J.~M. Zain, H.~Z. Mojahid, A.~K. Jumaat, N.~Hamdan, and D.~K. Bagal.
\newblock Calculating {SIRD} model parameters and re-susceptibility in
  {Malaysia}, {Pakistan}, {India} and {United} {Arab} {Emirates}.
\newblock {\em Communications in Nonlinear Science and Numerical Simulation},
  2025.
\newblock 108902.

\bibitem{Verma2021}
M.~Verma, A.~K. Verma, and A.~K. Misra.
\newblock Mathematical modeling and optimal control of carbon dioxide emissions
  from energy sector.
\newblock {\em Environment, Development and Sustainability},
  23(9):13919–13944, 2021.

\bibitem{Korda2017}
M~Korda and I.~Mezić.
\newblock On convergence of extended dynamic mode decomposition to the
  {K}oopman operator.
\newblock {\em Journal of Nonlinear Science}, 28(2):687–710, 2017.

\bibitem{Jin2024}
Y.~Jin, L.~Hou, and S.~Zhong.
\newblock Extended dynamic mode decomposition with invertible dictionary
  learning.
\newblock {\em Neural Networks}, 173:106177, 2024.

\bibitem{Mickens1993}
R.~E. Mickens.
\newblock {\em Nonstandard finite difference models of differential equations}.
\newblock World Scientific, 1993.

\bibitem{Ehrhardt2013}
M.~Ehrhardt and R.~E. Mickens.
\newblock A nonstandard finite difference scheme for convection–diffusion
  equations having constant coefficients.
\newblock {\em Applied Mathematics and Computation}, 219(12):6591--6604, 2013.

\bibitem{Costa2024}
G.~Costa, M.~Lobosco, M.~Ehrhardt, and R.~Reis.
\newblock Mathematical analysis and a nonstandard scheme for a model of the
  immune response against {COVID-19}.
\newblock {\em Mathematical and Computational Modeling of Phenomena Arising in
  Population Biology and Nonlinear Oscillations}, page 251–270, 2024.

\bibitem{Maamar2023}
M.~H. Maamar, M.~Ehrhardt, and L.~Tabharit.
\newblock A nonstandard finite difference scheme for a time-fractional model of
  {Zika} virus transmission.
\newblock {\em Mathematical Biosciences and Engineering}, 21(1):924–962,
  2023.

\bibitem{JHU2020}
Johns~Hopkins University.
\newblock {COVID-19} dashboard; \url{https://coronavirus.jhu.edu/data},
  Accessed on: August 18, 2025 (2020).

\bibitem{CDCCOVID}
Centers for Disease~Control and Prevention.
\newblock {CDC COVID} data tracker;
  \url{https://www.cdc.gov/covid-data-tracker/index.html}, Accessed on: August
  18, 2025.

\bibitem{CDCFlu}
Centers for Disease~Control and Prevention.
\newblock Seasonal influenza ({Flu}), \url{https://www.cdc.gov/flu/index.htm},
  Accessed on: August 18, 2025.

\bibitem{ECDCFlu}
European~Centre for Disease~Prevention and Control.
\newblock Seasonal influenza data;
  \url{https://www.ecdc.europa.eu/en/seasonal-influenza}, Accessed on: August
  18, 2025.

\bibitem{CDCEbola2014}
Centers for Disease~Control and Prevention.
\newblock Ebola virus disease outbreak, {West Africa}, october 2014;
  \url{https://www.cdc.gov/mmwr/preview/mmwrhtml/mm6343a3.htm}, Accessed on:
  August 18, 2025.

\bibitem{WHOEbola2014}
World~Health Organization.
\newblock Ebola outbreak 2014–2016, {West Africa};
  \url{https://www.who.int/emergencies/situations/ebola-outbreak-2014-2016-west-africa},
  Accessed on: August 18, 2025.

\bibitem{CDCMeasles}
Centers for Disease~Control and Prevention.
\newblock {Measles Surveillance Manual};
  \url{https://www.cdc.gov/surv-manual/php/table-of-contents/chapter-7-measles.html},
  Accessed on: August 18, 2025.

\bibitem{CDCPinkBook}
Centers for Disease~Control and Prevention.
\newblock {Pink Book – Measles};
  \url{https://www.cdc.gov/pinkbook/hcp/table-of-contents/chapter-13-measles.html},
  Accessed on: August 18, 2025.

\end{thebibliography}

\end{document}